\documentclass{amsart}
\usepackage{amsmath, amsthm, amssymb, amscd, mathrsfs, eucal, epsfig}

\begin{document}

\font\bfit=cmbxti10

\newtheorem{theorem}{Theorem}
\newtheorem{lemma}[theorem]{Lemma}
\newtheorem{sublemma}[theorem]{Sublemma}
\newtheorem{proposition}[theorem]{Proposition}
\newtheorem{corollary}[theorem]{Corollary}
\newtheorem{conjecture}[theorem]{Conjecture}
\newtheorem{question}[theorem]{Question}
\newtheorem{problem}[theorem]{Problem}
\newtheorem*{claim}{Claim}
\newtheorem*{criterion}{Criterion}
\newtheorem*{main_thm}{Main Theorem}

\theoremstyle{definition}
\newtheorem{definition}[theorem]{Definition}
\newtheorem{construction}[theorem]{Construction}
\newtheorem{notation}[theorem]{Notation}

\theoremstyle{remark}
\newtheorem{remark}[theorem]{Remark}
\newtheorem{example}[theorem]{Example}

\numberwithin{equation}{subsection}

\def\H{\mathbb H}
\def\Z{\mathbb Z}
\def\N{\mathbb N}
\def\R{\mathbb R}
\def\Q{\mathbb Q}
\def\D{\mathcal D}
\def\E{\mathcal E}
\def\RR{\mathcal R}
\def\RP{\mathbb{RP}}
\def\P{{\bf P}}
\def\F{\mathcal F}
\def\T{\mathcal T}
\def\ds{\displaystyle}
\def\l{\textnormal{l}}
\def\sl{\textnormal{sl}}

\def\ab{\textnormal{ab}}
\def\CAT{\textnormal{CAT}}
\def\cl{\textnormal{cl}}
\def\scl{\textnormal{scl}}
\def\homeo{\textnormal{Homeo}}
\def\rot{\textnormal{rot}}
\def\area{\textnormal{area}}
\def\inte{\textnormal{int}}		
\def\hom{\textnormal{Hom}}
\def\Aut{\textnormal{Aut}}

\def\Id{\textnormal{Id}}
\def\id{\textnormal{id}}
\def\SL{\textnormal{SL}}
\def\Sp{\textnormal{Sp}}
\def\PSL{\textnormal{PSL}}
\def\length{\textnormal{length}}
\def\fill{\textnormal{fill}}
\def\rank{\textnormal{rank}}
\def\til{\widetilde}

\title{Quasimorphisms and laws}
\author{Danny Calegari}
\address{Department of Mathematics \\ Caltech \\
Pasadena CA, 91125}
\email{dannyc@its.caltech.edu}

\date{1/11/2010, Version 0.02}

\begin{abstract}
Stable commutator length vanishes in any group that obeys a law.
\end{abstract}

\maketitle

If $G$ is a group and $g$ is an element of the commutator subgroup $[G,G]$, 
the {\em commutator length} of $g$, denoted $\cl(g)$, is the least number of commutators in $G$
whose product is $g$. The {\em stable commutator length}, denoted $\scl(g)$, 
is the limit $\scl(g):= \lim_{n \to \infty} \cl(g^n)/n$.

A group $G$ is said to obey a {\em law} if there is a free group $F$ (which may be assumed to have finite
rank) and a nontrivial element $w \in F$
so that for every homomorphism $\rho:F \to G$, we have $\rho(w)=\id$. For example, abelian (or, more
generally, nilpotent or solvable) groups obey laws. The free Burnside groups $B(m,n)$
with $m\ge 2$ generators and odd exponents $n \ge 665$ are perhaps the best known
examples of non-amenable groups that obey laws; see e.g. Adyan \cite{Adyan}.

The point of this note is to prove the following:

\begin{main_thm}
Let $G$ be a group that obeys a law. Then $\scl(g)=0$ for every $g \in [G,G]$.
\end{main_thm}

The proof is very short, given some basic facts about stable commutator length, which we recall for the
convenience of the reader. A basic reference is
Bavard's paper \cite{Bavard} or the monograph \cite{Calegari_scl}, especially Chapter~2.

\begin{definition}
A {\em homogeneous quasimorphism} on a group $G$ is a function $\phi:G \to \R$ that restricts to a homomorphism
on every cyclic subgroup, and for which there is a least number
$D(\phi)\ge 0$ (called the {\em defect}) so that for any $g,h \in G$ there is an inequality
$|\phi(gh)-\phi(g)-\phi(h)| \le D(\phi)$.
\end{definition}

The defect satisfies the following formula:
\begin{lemma}[\cite{Bavard}, Lem.~3.6 or \cite{Calegari_scl}, Lem.~2.24]\label{defect_lemma}
Let $\phi$ be a homogeneous quasimorphism. Then there is an equality
$$\sup_{g,h \in G} |\phi([g,h])| = D(\phi)$$
\end{lemma}

Bavard duality (see e.g. \cite{Bavard} or \cite{Calegari_scl}, Thm.~2.70) says that for any $g \in [G,G]$,
there is an equality $\scl(g) = \sup_\phi \phi(g)/2D(\phi)$ where the supremum is taken over all
homogeneous quasimorphisms $\phi$ with nonzero defect. In particular, $\scl$ is nontrivial on $G$
if and only if $G$ admits a homogeneous quasimorphism with nonzero defect. 

On the other hand, there is a topological formula for $\scl$.
Let $X$ be a space with $\pi_1(X)=G$, and let $\gamma:S^1 \to X$ be a free homotopy class representing
the conjugacy class of $g \in G$. If $\Sigma$ is a compact, oriented surface without sphere or disk
components, a map $f:\Sigma \to X$ is {\em admissible} if the map $\partial f:\partial \Sigma \to X$ can be
factorized as $\partial \Sigma \xrightarrow{d} S^1 \xrightarrow{\gamma} X$. For an admissible map,
define $n(\Sigma)$ by the equality $d_*[\partial \Sigma] = n(\Sigma)[S^1]$ in $H_1$; i.e. $n(\Sigma)$ is the degree with
which $\partial \Sigma$ wraps around $\gamma$. By reversing the orientation of $\Sigma$ if necessary, we assume
$n(\Sigma)\ge 0$. With this notation, one has the following formula:
\begin{lemma}[\cite{Calegari_scl}, Prop.~2.10]\label{chi_lemma}
With notation as above,
$$\scl(g) = \inf_\Sigma \frac {-\chi(\Sigma)} {2n(\Sigma)}$$
where $\chi$ denotes Euler characteristic, and the infimum is taken over all compact, oriented
surfaces and all admissible maps.
\end{lemma}
Notice that both $\chi(\cdot)$ and $n(\cdot)$ are multiplicative under finite covers.

We are now in a position to prove the main theorem.
\begin{proof}
Suppose that $G$ obeys a law. Then there is a free group $F$ and a nontrivial word $w \in F$ so that
any homomorphism from $F$ to $G$ sends $w$ to $\id$. Let $F_2$ be free on generators $x,y$. We can
embed $F$ in $F_2$, and express $w$ as a word $v$ in the generators $x,y$. Hence any homomorphism
from $F_2$ to $G$ sends $v$ to $\id$.

Let $X$ be a space with $\pi_1(X)=G$. Let $\Sigma$ be a once-punctured torus. We choose generators for
$\pi_1(\Sigma)$, and identify this group with $F_2 = \langle x,y\rangle$. Let $\alpha$ be a loop on
$\Sigma$ whose free homotopy class represents the conjugacy class of $v$. Then any continuous map
$f:\Sigma \to X$ sends $\alpha$ to a null-homotopic loop.

Now suppose contrary to the theorem that $\scl$ does not vanish on $[G,G]$.
By Bavard duality there is a homogeneous quasimorphism $\phi$ with nonzero defect. Scale $\phi$
to have $D(\phi)=1$. Then by Lemma~\ref{defect_lemma}, for any $\epsilon>0$ there are elements $g,h$
in $G$ with $\phi([g,h])\ge 1-\epsilon$, and consequently $\scl([g,h])\ge 1/2 - \epsilon/2$ by
Bavard duality. 

Let $\gamma:S^1 \to X$ be a loop representing the conjugacy class of $[g,h]$. There is a map $f:\Sigma \to X$
whose boundary represents the free homotopy class of $\gamma$. As above, the loop $\alpha$
on $\Sigma$ maps to a null-homotopic loop in $X$. By \cite{Scott}, there is a finite cover $\til{\Sigma}$
of $\Sigma$ of degree $d$ (depending on $\alpha$), so that some lift $\til{\alpha}$ of $\alpha$ is homotopic
to an embedded loop $\alpha'$. Composing the
covering map with $f$ gives a map $\til{f}:\til{\Sigma} \to X$ for which $\til{f}(\alpha')$ is
null-homotopic in $X$. Since $\alpha'$ is embedded, we can compress $\til{\Sigma}$ along $\alpha'$
to produce a new surface $\Sigma'$ mapping to $X$ by $f'$. The map $f'$ is admissible for $\gamma$, and
satisfies $n(\Sigma')=d$. Moreover, $\chi(\til{\Sigma}) = -d$, and $\chi(\Sigma') = 2-d$.
Consequently, by Lemma~\ref{chi_lemma}, we have $\scl([g,h]) \le 1/2 - 1/d$.

Since $d$ is fixed (depending only on the law satisfied by $G$) but $\epsilon$ is arbitrary, we
obtain a contradiction. Hence $\scl$ vanishes identically on $[G,G]$, as claimed.
\end{proof}

\begin{remark}
The statement of the main theorem may be rephrased positively as saying that if $\scl$ is nonzero
on $G$, then for any positive integer $n$, there are homomorphisms $F_2 \to G$ which are injective 
on the ball of radius $n$.
\end{remark}

If $w$ is a word in a free group $F$, define a {\em $w$-word} in $G$ to be the image of $w$ under
a homomorphism $F \to G$. Let $G(w)$ be the subgroup of $G$ generated by $w$-words. The
{\em $w$-length} of $g \in G(w)$, denoted $\l(g|w)$, 
is the smallest number of $w$-words and their inverses whose product is $g$ (commutator length is the case
$w = xyx^{-1}y^{-1} \in \langle x,y\rangle$),
and the {\em stable $w$-length}, denoted $\sl(g|w)$ is $\sl(g|w):=\lim_{n \to \infty} \l(g^n|w)/n$.
\begin{question}
Is there an example of a group that obeys a law, but for which $\sl(\cdot|w)$ is nontrivial for some
$w$?
\end{question}

This work was partially supported by NSF grant DMS 0707130. I would like to thank the anonymous
referee for helpful comments on an earlier draft.


\begin{thebibliography}{99}

\bibitem{Adyan}
  S. Adyan,
  \emph{Random walks on free periodic groups},
  Izv. Akad. Nauk SSSR Ser. Mat. {\bf 46} (1982), no. 6, 1139--1149
\bibitem{Bavard}
  C. Bavard,
  \emph{Longueur stable des commutateurs},
  Enseign. Math. (2) {\bf 37} (1991), no. 1-2, 109--150
\bibitem{Calegari_scl}
  D. Calegari,
  \emph{scl},
  MSJ Memoirs {\bf 20}. Mathematical Society of Japan, Tokyo, 2009
\bibitem{Scott}
  P. Scott,
  \emph{Subgroups of surface groups are almost geometric},
  J. London Math. Soc. (2) {\bf 17} (1978), no. 3, 555--565
\end{thebibliography}
\end{document}